\documentclass[12pt,a4paper]{article}

\usepackage{amssymb,latexsym}
\usepackage{amsmath}
\usepackage{fancyvrb}
\usepackage{graphicx}
\usepackage{algorithmic}
\usepackage{listings}
\usepackage{boxedminipage}
\usepackage{fancybox}
\usepackage{amsthm}
\usepackage{subfig}

\newcommand{\Bc}{\begin{center}}
\newcommand{\Ec}{\end{center}}

\theoremstyle{plain}
\newtheorem{Theorem}{Theorem}[section]

\title{A Barzilai-Borwein $l_1$-Regularized Least Squares Algorithm for Compressed Sensing}
\author{R. Broughton, I. Coope, P. Renaud, R. Tappenden*}
\date{}
\begin{document}
\maketitle
\begin{abstract}
Problems in signal processing and medical imaging often lead to calculating sparse solutions to under-determined linear systems. Methodologies for solving this problem are presented as background to the method used in this work where the problem is reformulated as an unconstrained convex optimization problem. The least squares approach is modified by an $l_1$-regularization term. A sparse solution is sought using a Barzilai-Borwein type projection algorithm with an adaptive step length. New insight into the choice of step length is provided through a study of the special structure of the underlying problem.  Numerical experiments are conducted and results given, comparing this algorithm with a number of other current algorithms.
\end{abstract}

\section{Introduction}

Many problems in signal processing and medical imaging can be described by the following linear model,
\begin{eqnarray}
\label{noisy}
b = Ax + v,
\end{eqnarray}
where $A \in \mathbb{R}^{m \times n}$ ($m < n$), $b \in \mathbb{R}^m$ is a vector of observations, $x \in \mathbb{R}^{n}$ is the vector of unknowns and $v$ is a noise vector usually assumed to be Gaussian. The aim is to determine a \emph{sparse} solution $x$. This is an ill-posed problem because $A$ is under-determined. In the over-determined case a standard approach is to solve for $x$ in a least-squares sense by minimizing $\|Ax-b\|_2^2$. However in the under-determined case least-squares regression leads to over-fit. Therefore a standard technique in statistical and signal processing problems is to incorporate a regularization term. As the solution vector $x$ is known to be sparse, early work, (see for example \cite{Muir73}), suggest regularization with an $l_1$ term (rather than Tikhonov (or $l_2$) regularization \cite{Hansen98}). This leads to the unconstrained convex optimization problem,
\begin{eqnarray}
\label{formula}
\min_x \|Ax-b\|_2^2 + \lambda \|x\|_1.
\end{eqnarray}
Here $\lambda > 0$ is a regularization parameter. The value of the scalar $\lambda$ is important, for example, if $\lambda$ is too large then the solution is the trivial one $x=0$, (see \cite{Boyd07}). The introduction of the $l_1$-regularization term significantly promotes a sparse solution while maintaining the convexity of the objective function.\\

In the next section we briefly mention some other approaches to finding sparse solutions to problem (\ref{noisy}) before focusing attention on specific implementations for solving problem (\ref{formula}). Section \ref{algorithm} introduces the highly successful variation on steepest descent proposed by Barzilai and Borwein \cite{Barzilai88} and this is exploited in the algorithm presented in section \ref{bbkey} with numerical results presented in section \ref{results}.\\

\section{Previous Approaches}
Several optimization algorithms have been recently proposed with the aim of determining a sparse $x$ satisfying (\ref{noisy}). Some notable approaches are discussed here.\\

In 2005 Cand$\grave{e}$s and Romberg \cite{Candes05} described an algorithm to solve the problem,
\begin{eqnarray}
\label{l1magic}
\min_x \|x\|_1 \quad s.t. \quad \|Ax-b\|_2^2 \leq \epsilon^2.
\end{eqnarray}
The use of the $l_1$-norm induces sparsity in $x$ while the constraint ensures $b \approx Ax$. (We recall that $b$ is observed in the presence of noise so it is reasonable not to enforce $b=Ax$ exactly).  The algorithm, so-called $l_1$-$magic$, is available online at \texttt{http://www.l1magic.org}.\\

More recently other groups have focused on devising algorithms for the solution of (\ref{formula}). A group at Stanford University \cite{Boyd07} began their work by formulating the dual. A new variable $z \in \mathbb{R}^{m}$ was introduced leading to the equivalent primal problem,
\begin{eqnarray*}
\min_{x,z} & \quad z^Tz + \lambda \|x\|_1\\
s.t. & \quad z = Ax - b.
\end{eqnarray*}
Dual variables $\nu_i$ were associated with the equality constraints $z_i$ and the Lagrange dual problem is
\begin{eqnarray*}
\max_{\nu} & \quad G(\nu) = -\frac{1}{4}\nu^T\nu - \nu^Tb\\
s.t.& \quad |(A^T\nu)_i| \leq \lambda.
\end{eqnarray*}
The primal problem (\ref{formula}) satisfies Slater's condition \cite{Boyd07} so the optimal value $f(x^*)$ of the primal problem is equal to that of the dual. The duality gap was used as a stopping criterion for their algorithm. (For more on convex duality see for example \cite{Boyd04}). Next (\ref{formula}) was transformed into the convex quadratic problem with linear inequality constraints:
\begin{eqnarray*}
\min_{x,u}& \quad \|Ax-b\|_2^2 + \lambda \sum_{i=1}^{n} u_i\\
s.t. & \quad -u_i \leq x_i \leq u_i \quad i=1,\dots,n.\frac{}{}
\end{eqnarray*}
An interior-point truncated Newton method was used to solve this problem. The Matlab code for this ($l_1$-$ls$) algorithm is publicly available online at\\
\texttt{http://www.stanford.edu/$\thicksim$boyd/l1\_ls/}.\\

A third group (Figueiredo, Nowak and Wright, \cite{Wright07}) reformulated (\ref{formula}) as the bound constrained quadratic programme,
\begin{eqnarray}
\label{BCQP1}
\min_{u,v}& \frac{1}{2}\|A(u-v) - b\|_2^2+\tau \sum_{i=1}^{n} u_i + \tau \sum_{i=1}^n v_i\\
s.t.& u,v \geq 0  \notag
\end{eqnarray}
where the substitution $x = u - v$, $u,v \geq 0$ has been made. A projected Barzilai-Borwein (PBB) \cite{Dai05} type method was used to determine an approximate solution of (\ref{BCQP1}). Matlab code for this (GPSR) algorithm is publicly available online at \texttt{http://www.lx.it.pt/$\thicksim$mtf/GPSR}.\\

Many other algorithms exist with applications to compressed sensing and the associated signal and image processing problems. For example, the SpaRSA (Sparse Reconstruction by Separable Approximation) algorithm \cite{Figueiredo07}, and the FISTA (Fast Iterative Shrinkage-Thresholding) algorithm \cite{Beck09}, are two very recent algorithms which are further considered in section \ref{results}. Other current algorithms include: a projected Barzilai-Borwein type algorithm with applications in computed tomography \cite{Wang07}; the Gradient Projection, GP, algorithm (and the Steplength Selection for Gradient Projection, GPSS, variant) \cite{Loris09}; a gradient descent algorithm which uses a thresholding step to encourage sparsity \cite{Garg09}; and an algorithm for a non-convex compressed sensing problem, \cite{Chartrand08}.

\subsection{A Proposed Approach}

The $l_1$-magic algorithm \cite{Candes05} for finding a sparse solution to problem (\ref{noisy}) has three levels of iteration (nested loops) and as a consequence, runs relatively slowly. When the problem is reformulated as (\ref{formula}), the algorithm in \cite{Boyd07} uses two levels of iteration while the approach in \cite{Wright07} uses only one level as they do not use a backtracking line search. \\

The approach proposed here also aims to determine a sparse
solution $x$ using problem formulation (\ref{formula}). A
Barzilai-Borwein type algorithm with an alternating
step-length, $\alpha$, is employed. This approach (known as
the Projected Alternating Barzilai-Borwein or PABB algorithm) is based  on recent work by Dai and Fletcher \cite{Dai05}
who have investigated a variant of the PABB method. They claim
that this implementation performs better than the PBB method in practice. \\

Our approach uses two levels of iteration, an outer loop defining a search direction and new candidate point $x$, and an inner backtracking line-search loop.
However, the backtracking line-search is included only as a safeguard,
(as suggested in \cite{Dai05} to prevent iterates cycling).
This algorithm only enters the back-tracking line-search loop
under certain conditions which in practice rarely arise.\\

Also, as in the case of the $l_1$-\emph{ls} and GPSR algorithms, our approach only requires matrix-vector products
involving $A$ and $A^T$. At each iteration there are only two
matrix-vector products --- one vector multiplication with $A$
and one with $A^T$ --- unless the inner loop is required in
which case there is an additional multiplication with $A$ in the
backtracking line search.  The computational effort is therefore kept
low in each iteration.\\

\section{A Reformulation of the Problem}
\label{algorithm}

By making the substitution $x = u-v$, problem (\ref{formula}) can be recast as the bound constrained quadratic programme,
\begin{eqnarray}
\label{BCQP}
\min_{u,v}& \|A(u-v) - b\|_2^2 + \lambda \sum_{i=1}^{n} u_i + \lambda \sum_{i=1}^n v_i\\
s.t.& u,v \geq 0. \notag
\end{eqnarray}
As (\ref{BCQP}) is now a differentiable problem, the associated gradient is,
\begin{eqnarray}
\label{gradient}
g = \begin{bmatrix} 2A^T(A(u-v) - b) + \lambda \mathbf{1} \\ -2A^T(A(u-v) - b) + \lambda \mathbf{1} \end{bmatrix}
\end{eqnarray}
(where $\mathbf{1}$ is a vector of ones) and the associated Hessian is
\begin{eqnarray}
\label{hessian}
H = 2\begin{bmatrix} A^TA  & -A^TA \\ -A^TA & A^TA \end{bmatrix}.
\end{eqnarray}
At the solution of problem (\ref{BCQP}) we have either $u_i=0$ or $v_i=0$. Problems (\ref{formula}) and (\ref{BCQP}), although different, share a common minimizer. We prefer to solve problem (\ref{BCQP}) as the objective function is now differentiable.\\

Another point to note, (as mentioned in \cite{Wright07}), is that the introduction of a shift, $u \leftarrow u + \Delta$ and $v \leftarrow v + \Delta$ leaves $x$ unchanged. The gradient (\ref{gradient}) is also independent of this shift although the objective function value (\ref{BCQP}) increases by $2 \lambda \Delta$. Therefore in the algorithm presented in section \ref{algorithm}, the value of the primal objective function is calculated using formula (\ref{formula}) rather than (\ref{BCQP}) as this gives a lower value of the objective function.\\

The Lagrange dual of primal problem (\ref{formula}) is
\begin{eqnarray}
\label{dual}
\max_{\nu} &\quad G(\nu) = - \frac{1}{4}\nu^T\nu - \nu^T b\\
s.t. &\quad |(A^T\nu)_i| \leq \lambda \notag
\end{eqnarray}
where $G(\nu)$ is the dual objective function. (This is derived in more detail in \cite{Boyd07}. For more on duality see for example \cite{Boyd04}, \cite{Fletcher91} or \cite{Nocedal06}). A dual feasible point $\nu$ gives a lower bound on the optimal value of the primal problem and therefore an indication of the error in the computed solution. Furthermore, as (\ref{formula}) satisfies Slater's condition, the optimal value of the primal problem is equal to the optimal value of the dual. Thus we define the duality gap to be
\begin{eqnarray}
\label{dgap}
\eta = \|Ax-b\|_2 + \lambda\|x\|_1 - G(\nu).
\end{eqnarray}
This can be used as a stopping criterion which is described later.

\section{Barzilai-Borwein Key Features}
\label{bbkey}

In 1988 Barzilai and Borwein devised a novel gradient method for optimization problems, \cite{Barzilai88}. This Barzilai-Borwein algorithm has the unusual property that at some iterates the function value increases. Despite this property, the algorithm performs very well in practice. In fact, forcing a monotonic decrease in function value at each iteration can seriously impair the practical performance of the algorithm, (see \cite{Dai05}).\\

There has also been much interest in this algorithm more recently: the implementation of dynamical retards \cite{Luengo03}, analysis of convergence properties \cite{Dai02} and the introduction of a cyclic Barzilai-Borwein variant \cite{Dai06}, (see also the review by Fletcher \cite{Fletcher01}). In the following subsections we introduce and discuss some of the key features of the PABB algorithm.

\subsection{Step Length and the Projection Operator}

Consider first the unconstrained case. One of the key points of the Barzilai-Borwein algorithm is the step length $\alpha$. The quasi-Newton equation is,
\begin{eqnarray}
\label{qN}
y_k = Hs_k,
\end{eqnarray}
where $y_k = g(x_k) - g(x_{k-1})$, $s_k = x_k - x_{k-1}$ and $H$ is the Hessian ($H = \nabla^2f(x)$). Suppose we approximate $H$ by the matrix $\alpha^{-1}I$ where $\alpha >0$. Solving $$\min_{\alpha} \|y_{k} - \alpha^{-1}s_{k}\|^2$$ gives
\begin{eqnarray}
\label{alpha1}
\alpha^{BB_1}_k = \frac{s^T_{k}s_{k}}{s^T_{k}y_{k}}.
\end{eqnarray}
Similarly using $\alpha I$ to approximate $H^{-1}$ and solving $$\min_{\alpha} \|\alpha y_{k} - s_{k}\|^2$$ yields
\begin{eqnarray}
\label{alpha2}
\alpha^{BB_2}_k = \frac{s^T_{k}y_{k}}{y^T_{k}y_{k}}.
\end{eqnarray}
Equations (\ref{alpha1}) and (\ref{alpha2}) give the two step lengths used in the Barzilai-Borwein algorithm. In the case of the problem expressed by (\ref{BCQP}) we have the following result.
\begin{Theorem}
For the function (\ref{BCQP}), if $A \in \mathbb{R}^{m \times n}$, ($m<n$), has orthonormal rows, then the Barzilai-Borwein step length (\ref{alpha2}) satisfies
$$\alpha^{BB_2}_k = \frac{1}{4}.$$
\begin{proof}
Let $x_k = u_k - v_k$ so that $\tilde{u} = u_k - u_{k-1}$ and let $\tilde{v} = v_k - v_{k-1}$. Then using (\ref{gradient}),
\begin{eqnarray*}
y_k = g(x_k) - g(x_{k-1}) = \begin{bmatrix} 2A^TA(\tilde{u}-\tilde{v}) \\ -2A^TA(\tilde{u}-\tilde{v}), \end{bmatrix}
\end{eqnarray*}
so
\begin{eqnarray*}
y_k^Ty_k &=& 8(\tilde{u}-\tilde{v})^TA^TAA^TA(\tilde{u}-\tilde{v})\\
&=& 8(\tilde{u}-\tilde{v})^TA^TA(\tilde{u}-\tilde{v})
\end{eqnarray*}
as $AA^T = I$. Using the quasi-Newton condition (\ref{qN}) we know that $s_k^Ty_k = s_k^THs_k$ where $H$ is the Hessian matrix (\ref{hessian}). So
\begin{eqnarray*}
s_k^THs_k &=& \begin{bmatrix} \tilde{u}^T & \tilde{v}^T \end{bmatrix}\begin{bmatrix} 2A^TA  & -2A^TA \\ -2A^TA & 2A^TA \end{bmatrix} \begin{bmatrix} \tilde{u} \\ \tilde{v} \end{bmatrix}\\
&=& 2(\tilde{u}-\tilde{v})^TA^TA(\tilde{u}-\tilde{v})
\end{eqnarray*}
Thus
\begin{eqnarray*}
\alpha^{BB_2}_k = \frac{s^T_{k}y_{k}}{y^T_{k}y_{k}}
=\frac{2(\tilde{u}-\tilde{v})^TA^TA(\tilde{u} -\tilde{v})}{8(\tilde{u} -\tilde{v})^TA^TA(\tilde{u}-\tilde{v})}
= \frac{1}{4}
\end{eqnarray*}
\end{proof}
\end{Theorem}

Now we return to the constrained optimization case. A second key feature of this Barzilai-Borwein variant is the projection operator. Because we have a constrained optimization problem (\ref{BCQP}), once a search direction and step have been determined the projection operator (defined below) ensures the new candidate point $x$ is feasible. If we define the feasible set of (\ref{BCQP}) to be
\begin{eqnarray*}
\Omega = \{x : lb \leq x \leq ub\}
\end{eqnarray*}
where $lb$ and $ub$ are lower and upper bounds respectively, then the projection operator onto $\Omega$ is
\begin{eqnarray}
\label{projop}
P_{\Omega}(x) = mid(lb,x,ub)
\end{eqnarray}
where $mid(lb,x,ub)$ is the vector whose $ith$ component is the median of the set $\{lb_i,x_i,ub_i\}$. This operator ensures any $x$ is kept within the feasible region.

\subsection{An Adaptive Non-monotone Line-Search}
\label{ls}
The algorithm we propose includes a backtracking line-search loop. This line-search was proposed by Dai and Fletcher \cite{Dai05} who commented, ``the method again has a reference function value $f_r$ and each iteration must improve on the reference value. The method involves a small integer parameter $L>0$, and $f_r$ is reduced if the method fails to improve on the previous best value of $f$ in at most $L$ iterations. We dispense with the requirement
\begin{eqnarray*}
f(x_k + \rho d_k) \leq f_r + \theta \rho g_k^T d_k
\end{eqnarray*}
(where $d_k$ is the search direction, $\rho >0$ is a decreasing sequence of values and $\theta \in (0,1)$), to obtain a \emph{sufficient reduction} in $f$, since in real computation any reduction is bounded uniformly away from zero by a small amount \dots and this is sufficient to ensure global convergence. We refer to this kind of line search as an \emph{adaptive non-monotone line search}.'' The update strategy is clarified by the following pseudo-code where initially $f_r = \infty$, and $f_c = f_{best} = f(x_1)$.
\Bc
\fbox{\begin{minipage}{8cm}
\begin{align*}
if \quad & f(x_k) < f_{best}\\
& f_{best} = f(x_k), \: f_c = f(x_k), \: l=0\\
else\\
& f_c = \max\{f_c,f(x_k)\},\: l = l+1\\
& if \:\: l=L\\
& \qquad f_r = f_c, \: f_c=f(x_k), \: l=0 \\
& end\\
end\\
\end{align*}
\end{minipage}}
\Ec
This code reduces the reference function value $f_r$ to the candidate function value $f_c$ if $f_{best}$ has not been improved upon after $L$ iterations. This is enough to enforce convergence while still allowing non-monotone behaviour.\\

The choice of parameter $L$ is important. It represents the number of iterations allowed before a function decrease is enforced. For example, $L=1$ implies that the function value must be decreased at each iteration (a monotonic decrease in the objective function). As mentioned in section \ref{bbkey}, forcing a decrease in the objective function can impair the practical performance of the Barzilai-Borwein algorithm. It is suggested in \cite{Dai05} that suitable choices are $L=4$ or $L=10$. Initial testing showed little difference between the two choices, and so in the numerical results presented in section \ref{results}, $L=4$ is used.\\

\subsection{Bounds on Allowable Step length}
\label{sbound}
In an optimization problem we would like to take the step
\begin{eqnarray}
\label{nextx}
x_{k+1} = x_k + \alpha p_k
\end{eqnarray}
where $\alpha$  is the step length (given by either (\ref{alpha1}) or (\ref{alpha2})), and $p_k$ is the search direction. The objective function (\ref{BCQP}) is not strictly convex, that is, has a positive semi-definite Hessian. The original Barzilai-Borwein convergence theory applied to strictly convex quadratic functions and therefore extra safeguards on $\alpha$ may be needed to account for zero curvature.\\

In the strictly convex, quadratic, unconstrained case, the step lengths are automatically bounded by the reciprocal of the smallest eigenvalue of the Hessian, \cite{Raydan93}. In the present context the reciprocal of the smallest eigenvalue leads to an infinitely large step length so \cite{Raydan97} discuss the use of an upper bound $\alpha_{max} = 10^{30}$ as a safeguard which is implemented if the algorithm finds a direction of zero curvature, (or near zero curvature) at any iterate. However $\alpha_{max} = 10^{30}$ is not desirable in the current context when the solution is known to satisfy  an a priori bound, even in the presence of zero curvature. Hence we propose an upper bound on all iterates $x_k$ as follows.\\

A more appropriate crude upper bound for the solution vector is outlined here. At the unique minimizer $x^{*}$ we have,
\begin{eqnarray*}
f(x^{*}) \leq f(x) = \|Ax - b\|_2^2 + \lambda \|x\|_1
\end{eqnarray*}
for any $x$. Putting $x=0$ gives,
$$\lambda \|x^* \|_1 \leq  f(x^*) \leq b^Tb,$$
so that
$$\|x^*\|_1 \leq \frac{b^Tb}{\lambda}$$
and therefore
\begin{eqnarray}
\label{ub}
|x_{i}^*| \leq \frac{b^Tb}{\lambda}.
\end{eqnarray}

This is an \textit{a priori} bound on the components of $x^*$ and hence an \textit{a priori} upper bound on each $u_i$ and $v_i$. It is possible to improve this bound dynamically but numerical trials suggest this is not worthwhile. In any case we do not expect this bound to be active at the solution.\\

Another result discussed below also supports the inclusion of an upper bound as a safeguard against overly large step lengths. Let $B = A^TA$. Then the Hessian matrix (\ref{hessian}) can be written as the Kronecker product
\begin{eqnarray*}
H =  B \otimes \begin{bmatrix} 2 & -2\\ -2 & 2 \end{bmatrix}
\end{eqnarray*}
Using the properties of the eigenvalues of a Kronecker product (see \cite{Bernstein05}) $H$ has $2n-m$ zero eigenvalues corresponding to $2n-m$ directions of zero curvature. The remaining $m$ positive eigenvalues of $H$ are given by the positive eigenvalues of $4B$. In the simple case when A has orthonormal rows, the positive eigenvalues of $B$ are all equal to 1.\\

The dimension of the subspace of directions of zero-curvature is high. Thus there is an increased probability of encountering search directions for which the change in gradient would be tiny resulting in very large values of $\alpha$ for the next iteration. This supports the inclusion of the upper bound (\ref{ub}) which the solution is known to satisfy, which seems more appropriate than the $\alpha_{max} = 10^{30}$ approach.\\

\subsection{A Stopping Criterion}

An appropriate stopping criterion for any optimization algorithm is paramount to ensure that an accurate solution is located. A standard approach is to terminate when the norm of the projected gradient (see for example, \cite{Dai05}), is sufficiently small, indicating a stationary point has been found. The approach we favour is to use the duality gap as an indication of distance from the correct solution. So our termination criterion is
\begin{eqnarray}
\label{reldgap}
\frac{\eta}{G(\nu)} < tol
\end{eqnarray}
where $G(\nu)$ and $\eta$ are defined in (\ref{dual}) and (\ref{dgap}) respectively and $tol$ is some user-defined tolerance.\\

The stopping criterion (\ref{reldgap}) and a tolerance, $tol = 10^{-6}$, have been implemented in the BBCS algorithm outlined in this work and is used in all numerical results.

\subsection{A Barzilai-Borwein Algorithm for Compressed Sensing}
\label{algorithm}
Here we propose an algorithm based upon the ideas in the previous sections which aims to solve problem (\ref{BCQP1}). We refer to this algorithm as the Barzilai-Borwein algorithm for Compressed Sensing - BBCS algorithm. The BBCS algorithm is based on the algorithm described by Dai and Fletcher in \cite{Dai05} but has been tailored to problem (\ref{BCQP}) with tighter bounds on the allowable candidate vectors.\\

Recall the substitution $x = u-v$. Let
\begin{eqnarray*}
z = \begin{pmatrix} u\\ v\end{pmatrix}
\end{eqnarray*}
where $u$ and $v$ are defined as follows,
\begin{eqnarray*}
u = \max(x,0), \qquad v = \min(x,0)
\end{eqnarray*}
The steepest descent search direction is
\begin{eqnarray}
\label{sdx}
\hat{z} = z_k - \alpha_{k-1} g_k
\end{eqnarray}
where $g_k$ is the gradient defined by (\ref{gradient}) at the point $z_k$ and $\alpha$ is defined by either formula (\ref{alpha1}) or (\ref{alpha2}). As $\hat{z}$ may now violate the constraints, the projection operator (\ref{projop}) is used to give a point $z_{P}$ say, which is now feasible. The projection operator uses the upper bound $ub = b^Tb/\lambda$ as defined in section (\ref{sbound}) and a lower bound, $lb=0$. Thus the search direction $p$ used in the algorithm proposed here is
\begin{eqnarray*}
p = z_P - z_k.
\end{eqnarray*}
Based upon the previous arguments a backtracking line search loop may be used to encourage the algorithm to converge. That is, rather than forcing a monotonic decrease in function value at each iteration, a backtracking line search loop is used if the lowest function value $f_{best}$ has not been improved upon in the previous $L$ iterations. (Recall the discussion in section (\ref{ls})). The backtracking line-search is described by
\begin{eqnarray*}
z^+ = z_k + \beta p
\end{eqnarray*}
where $\beta = \frac{1}{2}, \frac{1}{4},\frac{1}{8},\dots$ until $f(z^+) < f_{r}$. Enter the adaptive non-monotone line search stage and update $f_r$, $f_c$ and $f_{best}$ according to the pseudo-code described in section (\ref{ls}). Finally,
\begin{eqnarray*}
x = u - v
\end{eqnarray*}
is computed and the duality gap (\ref{reldgap}) is monitored to check for convergence.\\

The results of this section are summarized in algorithmic form.\\

\begin{center}
\begin{minipage}{14cm}
\textbf{Step 0 (Initialization):} Set $\beta=1/2$, $L =4$, function reference values $f_{r}$, $f_c$ and $f_{best}$, $lb=0$, $ub= b^Tb/\lambda$ and $tol=10^{-6}$.\\
\textbf{Step 1:} Compute the step length $\alpha$ and gradient $g$.\\
\textbf{Step 2:} Compute $z$ and replace it with its projection, $mid(z-\alpha g, lb,ub)$.\\
\textbf{Step 3:} Compute the new search direction and $z_{k+1}$.\\
\textbf{Step 4:} If required, perform backtracking line-search and update the reference values: $f_r$, $f_c$ and $f_{best}$.\\
\textbf{Step 5:} Check whether the duality gap is sufficiently small. If so, terminate the algorithm, otherwise return to step 1.
\end{minipage}
\end{center}

\section{Numerical Results}
\label{results}

In this section we present numerical results obtained using the algorithm outlined in section \ref{algorithm}. These results illustrate the performance of the BBCS algorithm and how its performance compares with other recent algorithms -- namely, the $l_1$-$ls$ algorithm \cite{Boyd07}, the GPSR algorithm (both the monotone and the non-monotone versions) \cite{Wright07}, the SpaRSA algorithm (both the monotone and non-monotone versions) \cite{Figueiredo07}, and the FISTA algorithm \cite{Beck09}.\\

\subsection{A Sparse Signal Reconstruction Problem}
\label{spsignal}

\begin{figure}[where]
\label{Signal}
\centering
\includegraphics[width=8cm]{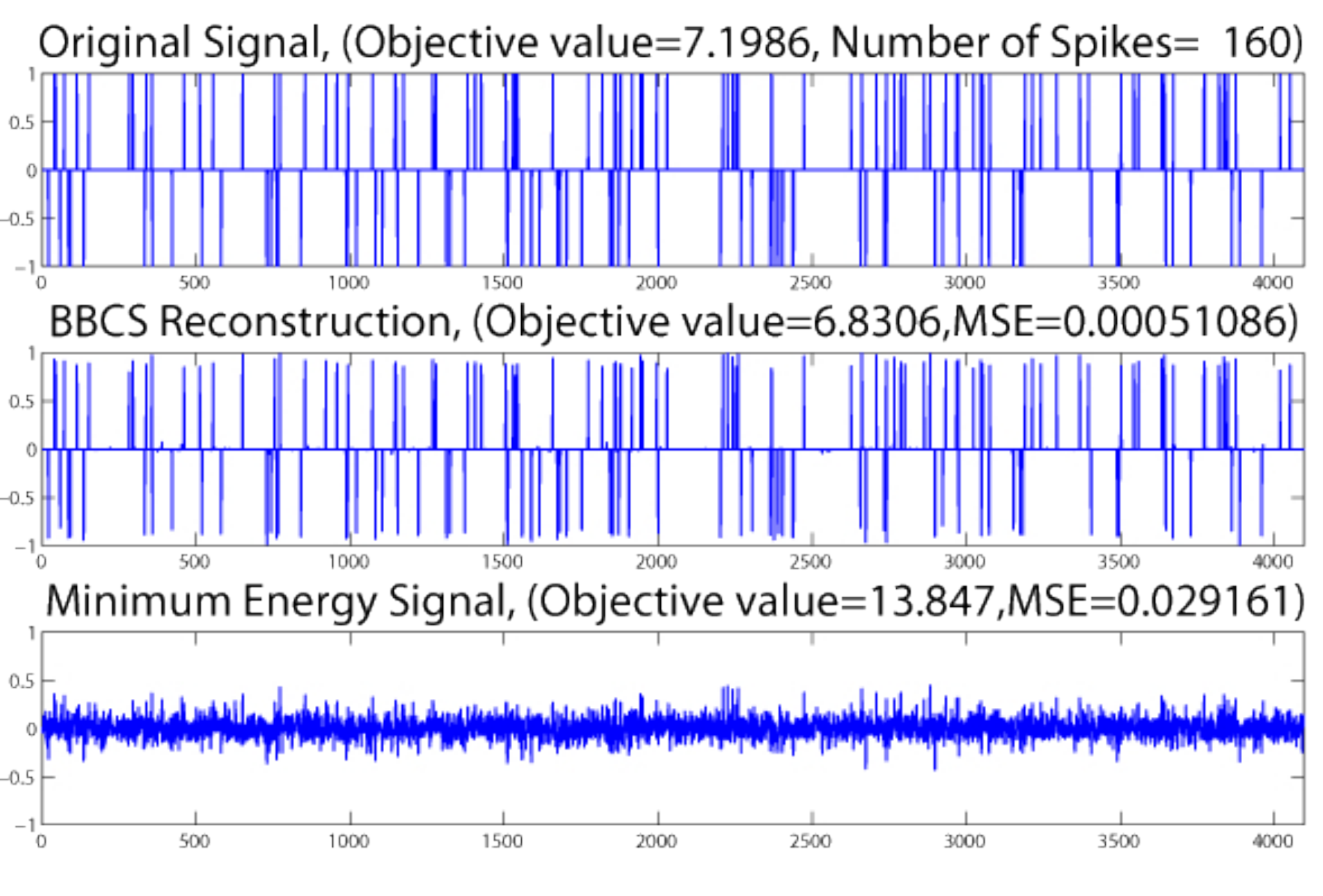}
\caption{Sparse signal reconstruction. From top to bottom: original signal, reconstruction from noisy observations, minimum energy solution.}
\end{figure}

The first numerical example demonstrated here is a sparse signal recovery experiment. A signal $x\in \mathbb{R}^{4096}$ consisting of 160 randomly placed spikes of amplitude $\pm 1$ was generated. A measurement matrix $A \in \mathbb{R}^{1024\times 4096}$ (representing $1024$ observations of the signal $x$) was constructed with Gaussian $\mathcal{N}(0,1)$ entries, and then the rows were orthonormalized (as for example in \cite{Candes05}, \cite{Boyd07}). The observation vector $b$ was formed according to (\ref{noisy}) where $v$ is drawn according to the Gaussian distribution with zero mean and variance $\sigma^2 = 10^{-4}$. The regularization parameter $\lambda$ was chosen to be
\begin{eqnarray*}
\lambda = 0.1\|A^Ty\|_{\infty},
\end{eqnarray*}
as this large $\lambda$ value seemed to encourage faster algorithm performance in initial numerical trials. As discussed in section \ref{ls}, the value $L=4$ was used along with a relative tolerance on the duality gap of $10^{-6}$. The step length $\alpha$ was computed using formula (\ref{alpha1}) except at every fourth iteration where formula (\ref{alpha2}) is used. The initial approximation was $x=\underline{0}$ where $\underline{0}$ is a vector of zeros. Figure (\ref{Signal}) shows the reconstruction results. The top plot shows the original signal. The middle plots shows the signal reconstructed using the BBCS algorithm. The BBCS algorithm  does an excellent job finding the positions of the non-zero components in the signal. The bottom plot shows the minimum energy solution (where $x= A^T(AA^T)^{-1}y$). Figure (\ref{Signal}) also shows the mean squared error, MSE, for both signal reconstructions\footnote{Here we follow \cite{Wright07} and define the MSE to be $\frac{1}{n}\|x_{true}-x\|_2^2$, where $x_{true}$ is the original signal.}. The signal reconstructed using the BBCS algorithm finds a solution with a low MSE indicating an accurate reconstruction.\\

Table (\ref{table1}) compares the runtimes of the Matlab implementation of our method and three existing methods on the problem described in subsection \ref{spsignal}. The BBCS algorithm is very efficient for this small problem.\\
\begin{table}[htp]
\label{table1}
\caption{CPU Times (Average Over 10 Runs) on the Experiment of Figure 1. (The subscript $m$ denotes the monotone version of the algorithm).}
\centering
\begin{tabular}{c c c }
\hline \hline
& & \\
Algorithm & & CPU Time (seconds)\\
& & \\
\hline
& &\\
BBCS & & 0.0770\\
BBCS$_m$ & & 0.0690\\
SpaRSA & & 0.0710\\
SpaRSA$_m$ & & 0.0810\\
GPSR & & 0.1230\\
GPSR$_m$ & & 0.0870\\
FISTA & & 0.1130\\
$l_1$-$ls$ & & 0.9080\\
& & \\
\hline
\end{tabular}
\end{table}

\subsection{Continuation}
Recent work in \cite{Wright07} and \cite{Figueiredo07} highlighted the possibility of implementing continuation schemes in their proposed algorithms. Here the algorithm starts with an initial regularization parameter $\lambda$ which is then reduced toward some desired value and the algorithm is warm-started for each successive value $\lambda$.\\

This scheme seems to have merit -- the algorithms with continuation schemes seem to find the solution to problem (\ref{BCQP}) faster then those without continuation schemes. However, we stress here that a continuation scheme could be applied to any algorithm as a means of improving speed. The results described here focus on \emph{the speed of the underlying algorithm} (while maintaining an accurate solution). Thus we compare the proposed algorithm without including a continuation scheme.\\

\subsection{Scalability Assessment}

An experiment proposed in \cite{Wright07} and \cite{Boyd07} aims to examine how the runtime of an algorithm changes as problem size grows. Their experiment is described as follows. Several random sparse matrices are considered whose entries are normally distributed. The dimensions of these matrices are $0.1n \times n$ where $n$ ranges from $10^4 - 10^6$. The sparsity of A is controlled to have $3n$ nonzero elements. For each data set, $x$ is also generated to be sparse with $n/4$ randomly placed components of length $\pm 1$. The measurements $Ax$ are corrupted with Gaussian noise of variance $\sigma^2 = 10^{-4}$. For each data set the regularization parameter is taken as $\lambda = 0.01\|A^Ty\|_{\infty}$.\\

\begin{figure}[where]
\label{SparseScale}
\centering
\includegraphics[width=8cm]{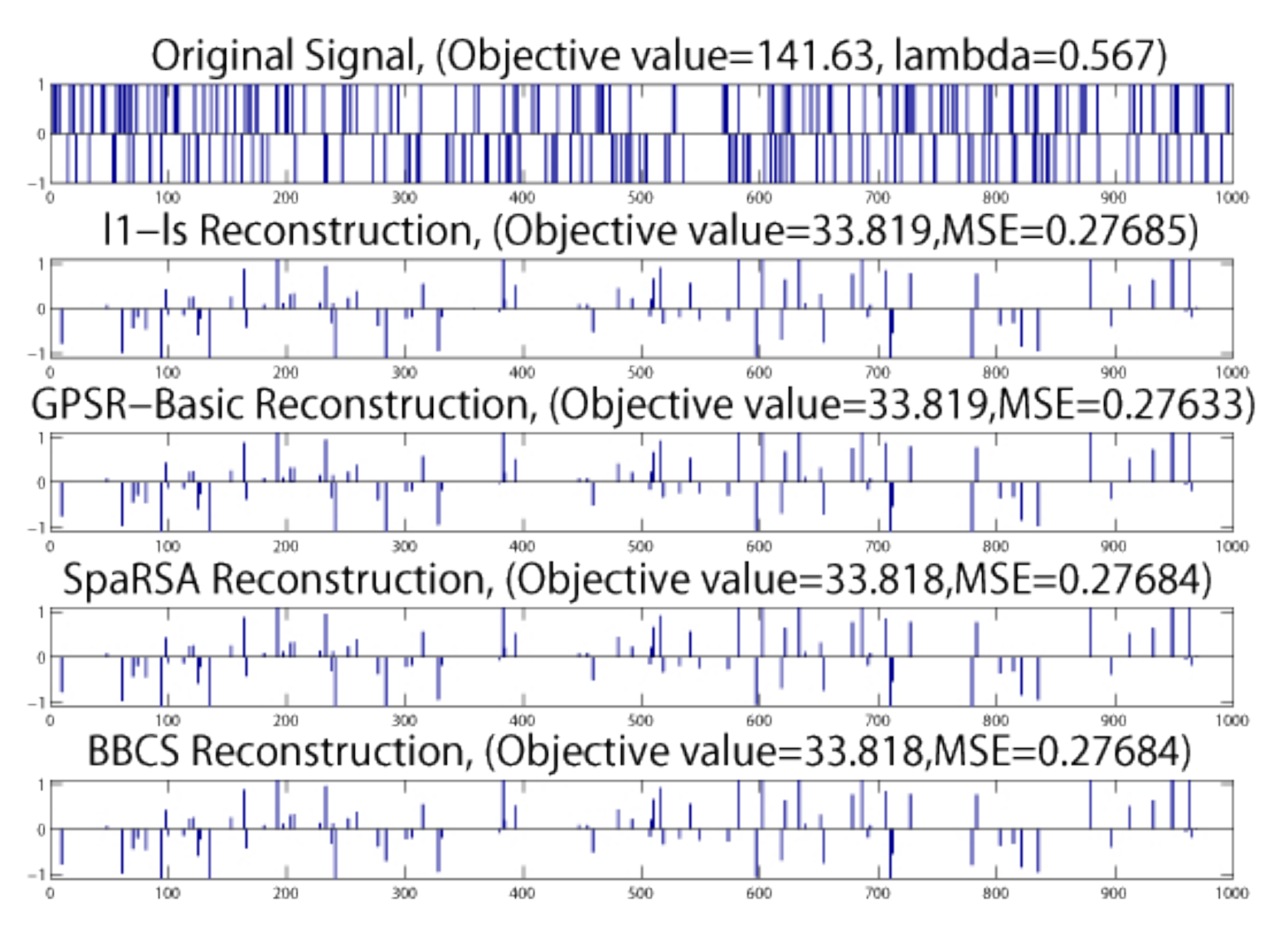}
\caption{A sparse matrix example. From top to bottom: Original signal, $l_1$-$ls$ signal reconstruction, GPSR-Basic reconstruction, SpaRSA reconstruction and the BBCS reconstruction. (The FISTA algorithm reconstructed $\hat{x} = \underline{0}$). Clearly none of the algorithms tested reproduce the original signal and all algorithms find an $\hat{x}$ vector giving a function value of $f(\hat{x}) = 33.8$ whereas the original signal gives a much higher function value of $f(x) = 141.6$.}
\end{figure}

An experiment based upon the above was implemented as a way of comparing the scalability of the GPSR, $l_1$-$ls$, FISTA, SpaRSA and BBCS algorithms. When performing this experiment we came across some interesting results. Figure (\ref{SparseScale}) shows a plot of the original signal and the signals reconstructed by the named algorithms on a problem of size $n=10^4$. It is clear from this figure that the algorithms are not reconstructing the original signal accurately. For the problem (\ref{formula}) each algorithm will always look for the sparsest solution and we know that we can always expect to find at least $n-m$ zeros in the solution. Because the problem is formulated with more than $m$ spikes, the algorithm chooses the solution vector which is sparsest and thus does not choose the original signal. We stress that these reconstructions are valid --- the algorithm is actually finding a solution vector $\hat{x}$ with $f(\hat{x}) << f(x)$ (where $x$ is the original signal), so from an optimization perspective the algorithms are working well. The problem is that the original signal is not being reconstructed. Since we are interested in reconstructing the original $x$ signal we have therefore decided to choose a scalability assessment based upon a non-random matrix.\\

The scalability assessment proposed here also considers the computational effort required as problem size increases. Observation matrices, (which are sub-matrices of a DCT matrix), were constructed. The dimensions of each matrix were $\frac{1}{8}n \times n$ with $n$ ranging from $2^{14} - 2^{20}$. Sparse signals with $\frac{n}{64}$ spikes of height $\pm 1$ were also generated for each matrix and the regularization parameter was chosen to be $\lambda = 0.1\|A^Ty\|_{\infty}$.\\

\begin{figure}[where]
\centering
\label{Scaleortho}
\includegraphics[width=8cm]{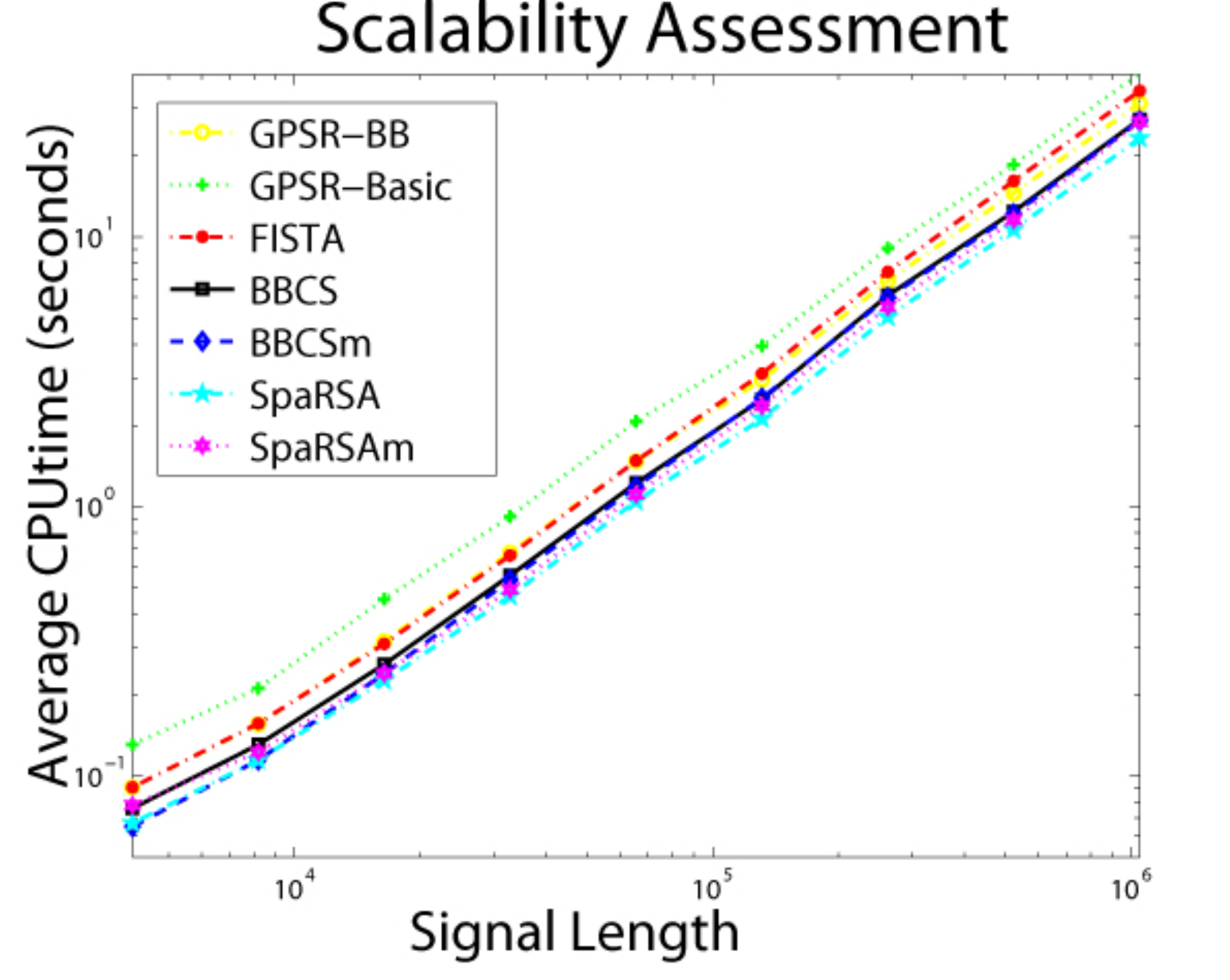}
\caption{Assessment the change in average CPU times for each algorithm as signal length increases.}
\end{figure}

All algorithms (BBCS, $l_1$-$ls$, GPSR, SpaRSA and FISTA), were tested using this experiment set-up. For each size $n$, ten matrices were generated and the average CPU time for each algorithm was found. The signal length vs average CPU times are shown in figure (\ref{Scaleortho}). From this we see that the BBCS algorithm is performing very competitively with the other algorithms.\\

The data from the scalability assessment was also used to estimate the computational complexity of each algorithm. That is, assume the computational cost is $\mathcal{O}(n^{\alpha})$ and estimate $\alpha$ based upon CPU times as $n$ increases. Table (\ref{table2}) gives the empirical estimates of the exponent $\alpha$. As shown, there is very little difference between the empirical computational complexity of each algorithm (around $1\%$). Table (\ref{table2}) also shows the CPU time for each algorithm when $n=2^{20}$. The BBCS (both monotone and non-monotone variants) are performing very competitively with the other algorithms.\\

\begin{table}[where]
\label{table2}
\caption{Empirical Estimate of the Exponent and Average CPU time of each algorithm (over 10 runs). Again the subscript $m$ denotes the monotone version of each algorithm.}
\centering
\begin{tabular}{c c c c c c}
\hline \hline
& & & & \\
Algorithm & & $\alpha$ Value & & CPU time ($n=2^{20}$)\\
& & & & \\
\hline
& & & & \\
BBCS & & 1.082 & & 27.32\\
BBCS$_m$ & & 1.108 & & 26.99\\
GPSR & & 1.053 & & 40.62\\
GPSR$_m$ & & 1.073 & & 31.35\\
SpaRSA & & 1.075 & & 23.29\\
SpaRSA$_m$ & & 1.078 & & 26.80\\
FISTA & & 1.096& & 34.85\\
& & & & \\
\hline
\end{tabular}
\end{table}

\section{Conclusion}

The problem of finding sparse solutions to large, under-determined linear systems in the presence of noise is an important one in signal processing, particularly in medical imaging. In this paper we have discussed a number of recent approaches and have proposed a variation of the PABB algorithm which we call the Barzilai Borwein algorithm for Compressed Sensing, BBCS, which provides safeguards in the case where the Hessian (\ref{hessian}) is positive semi-definite. These include the incorporation of an adaptive non-monotone line search, an upper bound on $x$ as a safeguard in the presence of zero curvature and a stopping criterion which provides a known bound on the error in our reconstruction.\\

The numerical results in Table (\ref{table1}) show that our algorithm is competitive with other existing algorithms. These results are encouraging because the underlying algorithm does not include any continuation schemes which would improve performance further.\\

As the scalability experiment shows, as the magnitude of the problem is increased our method retains accuracy and efficiency.\\

Future work includes the implementation of continuation schemes in the algorithm and investigating further the effects of signal reconstruction when the observation data matrix is sparse.

\section*{Acknowledgment}

Rachael Tappenden is the recipient of a Doctoral Scholarship from the New Zealand Institute of Mathematics and its Applications (NZIMA).

\providecommand{\bysame}{\leavevmode\hbox to3em{\hrulefill}\thinspace}
\providecommand{\MR}{\relax\ifhmode\unskip\space\fi MR }
\providecommand{\MRhref}[2]{%
  \href{http://www.ams.org/mathscinet-getitem?mr=#1}{#2}
}
\providecommand{\href}[2]{#2}

\end{document}